\documentclass[12pt]{amsart}
\usepackage{amscd,amsthm,amsmath} % SRC Specials: DVI [Inverse] Search
\newtheorem{thm}{Theorem}[section]
\newtheorem{cor}[thm]{Corollary}
\newtheorem{lem}[thm]{Lemma}
\newtheorem{prop}[thm]{Proposition}

\newtheorem{ques}[thm]{Question}

\theoremstyle{definition}
\newtheorem{defn}[thm]{Definition}
\theoremstyle{remark}
\newtheorem{rem}[thm]{Remark}
\numberwithin{equation}{section}
\begin{document}

\title{Equivariant Phantom maps}%
\author{Jianzhong Pan }%
\address{Institute of Math.,Academia Sinica ,Beijing 100080, China }%
\email{pjz@math03.math.ac.cn}%

\thanks{The author is partially supported by ZD9603 of Chinese Academy of Science
and  the NSFC projects 19701032 , 10071087  }%
\subjclass{55P91,55P60}%
\keywords{ Phantom map,equivariant homotopy}%

\date{Mar. 15, 2001}%
%\dedicatory{Ý}%
%\commby{Ý}%
% ----------------------------------------------------------------
\begin{abstract}
 A successful generalization of phantom map theory to the
 equivariant case for all compact Lie groups is obtained in this
 paper.
 One of the key observations is the discovery of the fact that
 homotopy fiber of equivariant completion splits as product of
 equivariant Eilenberg-Maclane spaces which seems impossible at
 first sight by the example of Triantafillou\cite{trian}.
\end{abstract}
\maketitle
% ----------------------------------------------------------------
 \section{Introduction} \label{S:intro}
A map $f:X \to Y$ from a CW complex $X$ is called a phantom map if
the restriction of $f$ to the skeleton $X_n$ , $f|X_n$ , is
trivial for each $n$. After the discovery of the first example of
phantom
 map by  Adams and Walker\cite{aw} , theory of phantom map receives a
  lot of attention . Among the many studies on phantom maps the most successful
  study is those conducted by Zabrodsky , McGibbon and Roitberg, see \cite{mcgib} for
  the references. Beside the study on the phantom map itself , remarkble
   applications of theory of phantom maps was given by
    Harper and Roitberg\cite{roitharp},\cite{roit2}  who applied it to
    compute $SNT(X)$ and $Aut(X)$,by Roitberg\cite{roit}
      and Pan\cite{pan} where several conjectures of McGibbon were settled, by
        Pan and Woo\cite{panwoo}  where a deep relation between  Tsukiyama problem about
self homotopy equivalence  and a generalization of phantom map was
found .

In the category of $G$-spaces with base point , a map $f:X \to Y$
from a $G$-CW complex $X$ is called a phantom map if the
restriction of $f$ to the skeleton $X_n$ , $f|X_n$ , is trivial
for each $n$ where $X_n$ is the $G$-n-skeleton $X$. Generalization
of theory of phantom maps to the equivariant case was first given
by Oda and Shitanda\cite{os1,os2} which is successful only when
the group is cyclic.  The trouble one meets when following the
standard approach to the theory of phantom map is the fact that
the argument depends strongly on the following well known result
in the nonequivariant homotopy theory.
\begin{thm} \label{T:splitt}
Let $X$ be a CW complex of finite type which is also a rational
H-space. Then $X$ is rationally equivalent to a product of
Eilenberg-Maclane spaces
\end{thm}
The corresponding result in equivariant case was obtained for
finite cyclic group by Triantafillou\cite{trian} and
counterexample for the above result in case the group is a product
of two cyclic groups in the same paper as cited above.

It is thus easy to understand why Oda and Shitanda's equivariant
generalization is successful only in case the group is cyclic.

A key observation of this paper is the following which makes it
possible the equivariant generalization of an alternative approach
to the theory of phantom maps.

\begin{thm} \label{T:key}
Let $Y$ be a $G$-CW complex of finite type and $c:Y \to \hat{Y}$
be the equivariant completion of $Y$ . Then $Y_{\rho}$ is homotopy
equivalent to a product of equivariant Eilenberg-Maclane spaces
where $Y_{\rho}$ is the homotopy fiber of $c$.
\end{thm}

 The organization of this paper is as follows.

 In section \ref{S:phanintro} definition of equivariant phantom
 map is recalled and several simple but needed later results are
 stated without proof since their proof is parallel to that of
 their nonequivariant analogue.

In section \ref{S:complet} we will study equivariant localization
and completion where the key Theorem\ref{T:key} among others will
be proved.

 In section\ref{S:phantom2}
The Theorem\ref{T:key} will be applied to compute the set of
homotopy classes of phantom maps from a $G$-CW complex $X$ to
another $G$-space $Y$ following a recent approach to the
computation in the nonequivariant case given by Roitberg and
Touhey\cite{roittou}. The paper ends with computation of
equivariant phantom maps between certain pairs of spaces.
 In concluding the Introduction , we 'd like to give the
following
\begin{ques}
Is the condition that  $Y$ is of finite type essential?
\end{ques}

\begin{rem}

1. The condition that  $Y$ is of finite type is indeed essential
to our argument in proving Theorem\ref{T:key}.

2. All the arguments in this paper generalizes directly to the
case of phantom pair as defined in \cite{os2}.
\end{rem}

In this paper, $G$ will denote a compact Lie group. For simplicity
we will assume throughout the paper that all $G$-spaces are based
$G$-1-connected although almost all results in this paper remains
true if we assume all $G$-spaces are nilpotent , the target is
$G$-simple and all fixed point set of the domain have finite
fundamental groups. $[X, Y]$ will be the set of based
nonequivariant homotopy classes of maps while $[X, Y]_G $ will be
he set of based equivariant homotopy classes of equivariant maps
when both $X$ and $Y$ are $G$-spaces. $map_*^G(X,Y)$ will be the
space of based equivariant maps between $G$-spaces $X$ and $Y$.

\section{Phantom maps:preliminaries}\label{S:phanintro}
Let's begin with  definition.
\begin{defn}
Let  $X$ be a $G$-CW  complex, $Y$ be a $G$-space  . Then a map
$f:X \to Y$ is called an $G$-phantom map if $f|X_{n}=0$ for all $n
\geq 0$. Denoted by
\[
Ph^G(X,Y)=\{f:X \to Y| f \text{ is an $G$-phantom map} \}
\]
\end{defn}

An equivariant analogue of Gray's description of $Ph(X,Y)$ is
true\cite{os2}
\[
Ph^G(X,Y)=\underset{\rightarrow}{\lim}^1_n[\Sigma X_n , Y]_G
\]
which follows from the following well known exact sequence
\[
* \to \underset{\rightarrow}{\lim}^1_n[\Sigma X_n , Y]_G \to [X,Y]_G \to
\underset{\rightarrow}{\lim}_n[ X_n , Y]_G \to *
\] \label{T:exact}

By the above formula we get the following according to Proposition
4.3 in \cite{mcgib} .
\begin{lem}\cite{os2} \label{T:complet}
Let $X$ and $Y$ be $G$-nilpotent $G$-CW complexes of finite type.
Then $$Ph^G(X,\hat{Y})=*$$
\end{lem}

Another result which we need is the following whose nonequivariant
analogue is well known.
\begin{prop} \label{T:ration}
Let $X$ be a $G$-CW complex and $Y$ be a  $G$-space of finite
 type. Then
 \[
 Ph^G(X,Y_{0})=*
 \]
\end{prop}
Actually this is an easy corollary of the following whose proof is
exactly the same

\begin{thm}\cite{stein}
Let $X$ be a $G$-CW complex and $Y$ be a  $G$-space of finite
 type. Let $\{X_{\alpha}\}$ be the set of all finite subcomplexes
 of $X$. Then

 $$[X,Y_{0}]_G \cong \underset{\leftarrow}{\lim}_{\alpha} [X_{\alpha} ,
 Y_{0}]_G$$
\end{thm}

\section{Equivariant localization and completion}\label{S:complet}

Let $X$ be a $G$-nilpotent . Then the equivariant localization and
completion are defined in \cite{may1},\cite{may2}. Denote by $l:X
\to X_{0} $ the equivariant rationalization and $c:X \to \hat{X}$
the equivariant completion.  Let $Y_{\rho}$ be the homotopy fiber
of $c$ and $X_{\tau}$ be the homotopy fiber of $l$. Then it is
easy to show that $X_{\tau} \to X \to X_{0}$ is a cofibration. Now
we can state first of our results concerning localization and
completion which is well known in nonequivariant case.
\begin{thm}\label{T:sull}
Let  $Y$ be $G$-CW complex of finite type. Then
\begin{itemize}
\item{$\pi_n(\hat{Y}^H)=\underset{p}{\prod} Ext(\mathbb{Z}_{p^{\infty}},
\pi_n(Y^H))$ for all closed subgroups $H$ of $G$.}
\item{For $W$ a finite CW complex, $c_*:[W,Y]_G \to [W,\hat{Y}]_G$ is injective}
\end{itemize}
\end{thm}
\begin{proof}
The first part can be found in Theorem 14 of \cite{may2}.

The second part can be proved by an argument similar to that of
Theorem 2.5.3 of \cite{hilton}.
\end{proof}
An immediate consequence of the above Theorem is
\begin{cor}\label{T:rho}
Let  $Y$ be a $G$-CW complex of finite type. Then $j:Y_{\rho} \to
Y $ is a $G$-phantom map
\end{cor}

\begin{prop}
Let $X$ be   $G$-CW complex   and $Y$ be $G$-space . Then the
followings hold:
\begin{itemize}
\item{$[\Sigma^n X_{(0)}, \hat{Y}]_G=*$ for all $n \geq 0$}
%,\tilde{H}^n(X_{(0)},\pi_i(\hat{Y}))=0 $  for all $n,i \geq 0$}
\item{$[\Sigma^n X_{\tau}, Y_{\rho}]_G=*$ for all $n \geq 0$}
% ,\tilde{H}^n(X_{\tau},\pi_i(Y_{\rho}))=0$ for all $n,i \geq 0$}
\end{itemize}
\end{prop}

For a proof see that of Proposition 2.3 in nonequivariant case in
\cite{panwoo2}.

\begin{prop}\label{T:4iso}
Let $X$ be  $G$-CW complex   and $Y$ be $G$-space. Then the
followings hold :
\begin{itemize}
\item{$\tau^*: [X,\hat{Y}]_G  \simeq [X_{\tau},\hat{Y}]_G$}
\item{$\rho_*:[X_{(0)},Y_{\rho}]_G  \simeq [X_{(0)},Y]_G$ }
\item{$c_*:[X_{\tau}, Y]_G  \simeq [X_{\tau}, \hat{Y}]_G$  }
\item{$l^*:[X_{(0)}, Y_{\rho}]_G  \simeq [X, Y_{\rho}]_G$  }
\end{itemize}
\end{prop}

\begin{proof}
First part of the theorem follows from Proposition 1,6 and Theorem
6 of \cite{may2} as noted by Oda and Shitanda\cite{os2}.

The second and the fourth parts follows from long homotopy exact
sequence.

For the proof of third part, it follows from the following exact
sequence that it suffices to prove the third part for each $X_n$
\[
* \to \underset{\rightarrow}{\lim}^1_n[\Sigma (X_n)_{\tau} , Y]_G \to [X_{\tau},Y]_G \to
\underset{\rightarrow}{\lim}_n[( X_n)_{\tau} , Y]_G \to *
\]
The above sequence is exact since an easy argument shows that
there is a cofibration for each $n$
\[
\bigvee G/H \times (S^n)_{\tau} \to ( X_n)_{\tau} \to (
X_{n+1})_{\tau}
\]
and $X_{\tau}$ is the direct limit of $( X_n)_{\tau}$.

An induction argument using above cofibration reduces the proof of
third part to that of a nonequivariant analogue of it which is of
course true.

\end{proof}

Now we proceed to prove Theorem\ref{T:key}. Let us record the
theorem again

\begin{thm}
%\label{T:key}
Let $Y$ be a $G$-CW complex of finite type and $c:Y \to \hat{Y}$
be the equivariant completion of $Y$ . Then $Y_{\rho}$ is homotopy
equivalent to a product of equivariant Eilenberg-Maclane spaces
where $Y_{\rho}$ is the homotopy fiber of $c$.
\end{thm}
\begin{proof}
To begin the proof ,  consider the following commutative diagram
following Roitberg and Touhey\cite{roittou}
\[
\begin{CD}
\Omega Y @>>\Omega c> \Omega \hat{Y} @>>d> Y_{\rho} @>> j>
Y @>>c> \hat{Y} \\
@VVV  @VVV  @VV\parallel V @VVrV @VVV \\
\Omega Y_{0} @>>> \Omega (\hat{Y})_{0} @>>d'> Y_{\rho} @>> j'>
Y_{0} @>>> (\hat{Y})_{0}
\end{CD}
\]
The right hand square is the pull-back by the $G$-Arithmetic
Square Theorem as in \cite{iwase} and horizontal sequences are
easy to see to be fibration sequences.

Now Corollary \ref{T:rho} implies that $j:Y_{\rho} \to Y$ is a
$G$-phantom and thus the composite $j'=l \circ j: Y_{\rho} \to
Y_{0}$  is $G$-phantom . Then Proposition \ref{T:ration} shows
that $j'\simeq *$ and thus $d':\Omega (\hat{Y})_{0} \to Y_{\rho}$
has a right $G$-homotopy inverse. It follows that $Y_{\rho}$ is an
equivariant H-space. As in nonequivariant case $Y_{\rho}$ is
rational. The proof completes if we proves that the equivariant
$k$-invaraints of $Y_{\rho}$  are all trivial.

 To prove it , note first that the above result is also true for
 $Y^{(n)}$ where  $Y^{(n)}$ is the equivariant $n$-th Postnikov section of
 $G$-space $Y$ \cite{elmen}. Thus ${(d')^{(n)}}: \Omega (\hat{Y^{(n)}})_{0} \to (Y^{(n)})_{\rho}$
has a right $G$-homotopy inverse.

  let's consider the following commutative diagram
 \[
 \begin{CD}
\Omega Y^{(n+2)} @>> \Omega c^{(n+2)}> \Omega \hat{Y}^{(n+2)} @>>d^{(n+2)}> (Y^{(n+2)})_{\rho} \\
@VVk_1V   @VVk_2V  @VVk_3V \\
K(\pi_1,n+2) @>>\Omega c^{(n+2)}_*> K(\pi_2,n+2) @>>d^{(n+2)}_*>
K(\pi_3,n+2)
 \end{CD}
 \]where $\pi_1=\pi_{n+1}(\Omega Y^{(n+2)})$ ,
$\pi_2=\pi_{n+1}(\Omega \hat{Y}^{(n+2)})$ and
$\pi_3=\pi_{n+1}((Y^{(n+2)})_{\rho})$. The fact that $\Omega
c^{(n)}$ is completion and the above commutative diagram imply
that $k_3 \circ d^{(n+2)}$ is trivial. Since
$\pi_{n+1}((Y^{(n+2)})_{\rho})$ is rational , it follows that $k_3
\circ {d'}^{(n+2)}$ is trivial. On the other hand ${d'}^{(n+2)}$
has a right $G$-homotopy inverse . It follows that $k_3$ is
trivial which completes the proof.

\end{proof}

\section{Phantom maps: localization and completion approach}\label{S:phantom2}
The most fundamental result in the theory of phantom map is the
following
\begin{thm}\label{T:fundamental}
Let $X$ be a $G$-CW  complex  of finite type , $Y$ be 1-connected
$G$-space of finite type. Then the followings are equivalent:
\begin{itemize}
\item{$f: X \to Y $ is a $G$-phantom map }
\item{the composite $X \to Y \to \hat{Y}$ is trivial}
\item{$f \simeq l \circ \tilde{f} where \tilde{f}: X_{(0)} \to Y $ }
\end{itemize}
\end{thm}

\begin{proof}
The proof follows easily from Proposition \ref{T:complet}, Theorem
\ref{T:sull} and \ref{T:4iso} as in \cite{os2}.
\end{proof}

Until now our argument follows closely the standard approach to
the theory of phantom map as developed by Zabrodsky\cite{za} and
this was what Oda and Shitanda would have been done. After that we
will face the computation problem. The standard approach suggests
to take $[X_{(0)}, Y]_G$ as an upper bound for $Ph^G(X,Y)$ and
compute it. Then the standard approach runs into trouble since the
equivariant analogue of Theorem \ref{T:splitt} is not true.  This
is the reason why   $Ph^G(X,Y)$ is computed only in case $G$ is a
finite cyclic group. Fortunately we have an alternative approach
to the computation first given by Roitberg and Touhey
\cite{roittou} . In their approach they used the completion
fibration. Now Theorem \ref{T:key} enables us to follow their
approach. Thus we have the following whose proof is parallel to
that of Theorem 1.1 of \cite{roittou} and will be omitted.

\begin{thm}\label{T:formul}
Let $X$ be a $G$-CW  complex  of finite type , $Y$ be 1-connected
$G$-space of finite type. Then
\[
Ph^G(X,Y) \approx [X, Y_{\rho}]_G/[X, \Omega \hat{Y}]_G
\]
 and

 \[
[X, Y_{\rho}]_G = \prod_{n \geq 1}H^n_G(X; \pi_{n+1}(Y) \bigotimes
\mathbb{R})
 \]
 where $\mathbb{R}$ is a vector space over $\mathbb{Q}$ of
 uncountable dimension.

 Furthermore if $X$ is a $G$-co-H-space or $Y$ is a $G$-H-space ,
 the above set theoretic bijections are group isomorphisms.
 Moreover the group structure on $Ph^G(X,Y)$ is abelian, divisible
 and independent of the $G$-co-H-structure on $X$  or the
 $G$-H-structure on $Y$.

\end{thm}

Now it's time to discuss when all maps between two spaces are
phantom.

As noted by  Oda and Shitanda \cite{os1}, The following Theorem is
especially useful in applying nonequivariant results to
equivariant case.
\begin{thm}\cite{elmen} \label{T:elmen}
Let $X$ be a $G$-CW complex and $Y$ be a $G$-space. If $[\Sigma^n
X^K , Y^H ]=*$ for all closed subgroups $K,H$ of $G$. Then
$[\Sigma^n X , Y]_G=*$ for all $n \geq 0$.
\end{thm}

As a corollary we have the following equivariant Miller-Zabrodsky
Theorem as noted by Oda and Shitanda \cite{os1}
\begin{thm}
Let $X$ be a $G$-CW complex with the following property:

There exists an integer $n$ such that, for all closed subgroups
$H$ of $G$, $\pi_i(X^H)=* $ for all $i \geq n$ and $\pi_i(X^H)$ is
locally finite for all $i < n$.

Then $map_*^G(X,\hat{Y})$ is weakly contractible  and thus
$$Ph^G(X,Y)=[X,Y]_G$$ for all finite
dimensional $G$-CW complex $Y$ .
\end{thm}

Combined with Theorem \ref{T:formul}, the above Theorem implies
the following
\begin{thm} \label{T:last}
Let $X$ be a $G$-CW complex of finite type and $Y$ be a finite
$G$-CW complex. If there exists an integer $n$ such that
$\pi_i(X^H)=* $ for all $i \geq n$ and all closed subgroups $H$ of
$G$, then
\[
Ph^G(X,Y)= \prod_{n \geq 1}H^n_G(X; \pi_{n+1}(Y) \bigotimes
\mathbb{R})
 \]
 where $\mathbb{R}$ is a vector space over $\mathbb{Q}$ of
 uncountable dimension.
\end{thm}

\begin{rem}
Equivariant Postnikov section of any $G$-CW complex satisfies the
codition for  $X$ in the above theorem . In particular any
Equivariant Eilenberg-MacLane space can be taken as $X$ in the
above theorem. This provides us with a large class of computable
phantom maps.

\end{rem}

Theorem \ref{T:elmen} combined with results in \cite{mcgib1}
yields another class of phantom maps
\begin{thm}
Let $K$ be a $G$-CW complex and $Y$ be a finite $G$-CW complex.
Let $X=\mathbf{Q}K$ where  $\mathbf{Q}K=\lim^n \Omega^n \Sigma^n
K$. Then $map_*^G(X,\hat{Y})$ is weakly contractible and thus
$$Ph^G(X,Y)=[X,Y]_G$$
\end{thm}

\begin{rem}
It follows from above theorem that the conclusion of Theorem
\ref{T:last} is also true.
\end{rem}

----------------------------------------------------------------

----------------------------------------------------------------


\begin{thebibliography}{99}
\bibitem{aw}
J.F.Adams and J.Walker,{\em An example in homotopy theory},Proc.
Camb. Phil. Soc., {\bf 60} (1960), pp. 699-700
%\bibitem{baues}
%H.J.Baues, {\em Rationale Homotopietypen}, Manus.Math. {\bf
%20}(1977), pp.119-131
%\bibitem{dror}
%E.Dror, W.G. Dwyer, D.M.Kan, {\em An arithmetic square for virtually nilpotent spaces}, Ill. J. Math., {\bf 21}(1977), pp.242-254
%\bibitem{dupont}
%N.Dupont {\em Problems and conjectures in rational homotopy
%theory}, Expos. Math., {\bf 12}(1994), pp. 323-352
\bibitem{elmen}
A.D. Elmendorf, {\em Systems of fixed point sets }, Trans. Amer.
Math. Soc. {\bf 277}(1983), pp.275--284
\bibitem{roitharp}
J.R.Harper, J.Roitberg, {\em Phantom maps and spaces of the same
n-type for all n }, J. of Pure and Applied Algebra  {\bf
80}(1992), pp.123--137
\bibitem{hilton}
P.Hilton, G.Mislin, J.Roitberg, {\em Localization of nilpotent groups and spaces}, North Holland Math. Series 15, 1975
\bibitem{iwase}
N.Iwase, {\em A continuous localization and completion}, Trans.
Amer. Math. Soc. 320(1990), no.1, pp.77-90


%\bibitem{kahn}
%D.W.Kahn, {\em Induced maps for Postnikov systems}, Trans.Amer.
%Math.Soc. {\bf 107}(1963), pp.432--450

\bibitem{may1}
J.P.May, J.McClure and G.Triantafillou, {\em Equivariant
localization}, Bull. London Math. Soc.  {\bf 14}(1982),
pp.223--230
\bibitem{may2}
J.P.May, J.McClure and G.Triantafillou, {\em Equivariant
completion}, Bull. London Math. Soc.  {\bf 14}(1982), pp.231--237
\bibitem{mcgib}
C.A.McGibbon,  {\em Phantom maps}, Handbook of Algebraic Topology,
North-Holland,1995, pp.1209--1258
\bibitem{mcgib1}
C.A.McGibbon,  {\em A note on Miller's theorem about maps out of
classifying spaces}, Proc. Am. Math. Soc.   {\bf 124}(1996),
pp.3241-3245
%\bibitem{miller}
%H.Miller, {\em The Sullivan conjecture on maps from classifying spaces}, Ann.Math., {\bf 120}(1984) , pp.39-87
\bibitem{os1}
N.Oda and Y.Shitanda,{\em Equivariant phantom Maps }, Publ. RIMS
kyoto Univ.  {\bf 24}(1988), pp.811-820

\bibitem{os2}
N.Oda and Y.Shitanda,{\em Localization,Completion and Detecting
Equivariant Maps on Skeletons }, Manuscript Math.,{\bf 65}(1989),
pp.1-18
\bibitem{pan}
Pan Jianzhong, {\em Having H-space structure is not a generic property},
submitted
\bibitem{panwoo}
Pan Jianzhong, Woo Mooha , {\em Phantom maps and injectivity of
forgetful maps, to appear in J.Japan. Math. Soc.}

\bibitem{panwoo2}
Pan Jianzhong, Woo Mooha , {\em Phantom elements and its
applications}, to appear in Contemp. Math. v.247 , 2001
\bibitem{roit2}
J.Roitberg, {\em Note on phantom phenomena and group of
self-homotopy equivalences}, Comment. Math. Helv. {\bf 66} (1991),
 pp. 448--457
 \bibitem{roit}
J.Roitberg, {\em The Lusternik-Schnirelmann category of certain
infinite CW complexes}, Topology {\bf 39}(2000), pp.95--101
\bibitem{roittou}
J.Roitberg, P. Touhey, {\em The homotopy fiber of profinite
completion}, Topology Appl. {\bf 103}(2000), pp.295--307
%\bibitem{schwartz}
%L.Schwartz,{\em Unstable modules over the Steenrod algebra and
%Sullivan fixed point set conjecture}, Chicago Lectures in
%Mathematics, 1994
\bibitem{stein}
R.J.Steiner, {\em Localization, completion and infinite
complexes}, Mathematika , {\bf 47}(1977), pp.1--15
\bibitem{trian}
G. Triantafillou, {\em Rationalization of Hopf $G$-spaces}, Math.
Z. {\bf 182}(1983), pp.485--500
\bibitem{za}
A.Zabrodsky, {\em On phantom maps and a theorem of H.Miller},
Israel J. Math. {\bf 58}(1987), pp.129--143



\end{thebibliography}
\end{document}